\documentclass[10pt,draft,twoside,leqno,a4paper]{article}
\usepackage{amsmath,amsthm,amsfonts,amssymb,verbatim}
\newtheorem{Theorem}{Theorem}[section]
\newtheorem{Lemma}{Lemma}[section]

\newtheorem{Proposition}{Proposition}[section]

\newcommand{\R}{\mathbb R}
\newcommand{\eps}{\varepsilon}

\newcommand{\dx}{\,{\mathrm d}x}
\newcommand{\dy}{\,{\mathrm d}y}
\newcommand{\weak}{\rightharpoondown}
\newcommand{\bra}{\langle}
\newcommand{\ket}{\rangle}

\newcommand{\us}{{\widetilde u}}

\newcommand{\Ge}{G_\eps}
\newcommand{\ue}{ u_\eps}

\newcommand{\Ie}{I_\eps}
\newcommand{\uo}{u_0}
\newcommand{\Io}{I_0}
\newcommand{\Ue}{U_\eps}
\newcommand{\Fe}{F_\eps}

\newcommand{\uj}{{u_j}}
\newcommand{\ujp}{u_j'}
\newcommand{\cj}{c_j}
\newcommand{\up}{\bar u}
\newcommand{\A}{\mathcal A}
\begin{document}
\title{Multiple vortices for a self-dual $CP(1)$ Maxwell-Chern-Simons model}
\author{ Francesco Chiacchio and Tonia Ricciardi\thanks{
Partially supported by 
the MIUR National Project {\em Variational Methods and
Nonlinear Differential Equations}} \\
{\small Dipartimento di Matematica e Applicazioni}\\
{\small Universit\`a di Napoli Federico II}\\
{\small Via Cintia, 80126 Naples, Italy}\\
{\small fax: +39 081 675665}\\
{\small \tt{francesco.chiacchio@dma.unina.it}}\\
{\small\tt{ tonia.ricciardi@unina.it}} }
\date{}
\maketitle
\begin{abstract}
We prove the existence of at least \emph{two} doubly periodic
vortex solutions for a
self-dual $CP(1)$ Maxwell-Chern-Simons model. 
To this end we analyze a system of two elliptic equations with exponential 
nonlinearities. Such a system is shown to be equivalent to a fourth-order 
elliptic equation admitting a variational structure.
\end{abstract}
\begin{description}
\item {\textsc{Key Words:}} nonlinear elliptic system, nonlinear
fourth-order elliptic equation, Chern-Simons vortex theory
\item {\textsc{MSC 2000 Subject Classification:}} 35J60
\end{description}
\setcounter{section}{-1}
\section{Introduction}
\label{sec:intro} 
The vortex solutions for the self-dual $CP(1)$ Maxwell-Chern-Simons model
introduced in \cite{KLL} (see also the monographs \cite{D,JT,Y}) are described by a system of two elliptic 
equations with exponential nonlinearities defined on a two-dimensional Riemannian manifold.
Such a system (henceforth, the ``$CP(1)$ system") was considered in \cite{CN}, where
among other results the authors prove 
the existence of one doubly periodic solution by super/sub methods.
On the other hand, formal arguments from physics as well as certain analogies with the $U(1)$
Maxwell-Chern-Simons model \cite{CK,RT} 
and with the $CP(1)$ ``pure" Chern-Simons model \cite{CN1,KLL1}
suggest that solutions to the $CP(1)$ system should be multiple.
In the special case of single-signed negative vortex points, 
a second solution for the $CP(1)$ system was exhibited in \cite{R3}. 
The method employed in \cite{R3} is not directly applicable
to the general case, due to the singularities produced
by the positive vortex points.
Our aim in this note is to prove 
multiplicity of solutions for the $CP(1)$ system in the general case of vortex points of \textit{either} sign.
In fact, we shall prove multiplicity for an abstract system of nonlinear elliptic equations 
which includes the $CP(1)$ system as a special case, 
thus emphasizing
some essential features of the $CP(1)$ system which ensure the multiplicity of vortex 
solutions. 
\par
For the sake of simplicity we define our equations 
on the flat 2-torus $M=\mathbb{R}^{2}/\mathbb{Z}^{2}$,
although it will be clear that corresponding results hold true on general compact Riemannian 
2-manifolds.  
We fix $p_1,\ldots,p_m\in M$ the ``positive vortex points" and $q_1,\ldots,q_n\in M$ 
the ``negative vortex points". The $CP(1)$ system as introduced in \cite{KLL} and analyzed in \cite{CN}
is given by
\begin{align*}
\Delta\us=&2q\left(-N+S-\frac{1-e^{\us}}{1+e^{\us}}\right)
-4\pi\sum_{j=1}^m\delta_{p_j}+4\pi\sum_{k=1}^{n}\delta _{q_k} &&\text{ on }M
\\
\Delta N=&-\kappa^2q^2\left(-N+S-\frac{1-e^{\us}}{1+e^{\us}}\right)
+q\frac{4e^\us}{(1+e^{\us})^2}N&&\text{ on }M
\end{align*}
where the couple $(\us,N)$ is the unknown variable, $q,\kappa>0$ and $S\in\R$ 
are constants and $\delta_{p_j}$,
$\delta_{q_k}$ are the Dirac measures centered at $p_j$, $q_k$.
Setting $v=N-S$, $s=-S$, $\lambda=2/\kappa$, $\eps=1/(\kappa q)$, the
above system takes the form:
\begin{align}
\label{sys_u} 
-\Delta\us=&\eps ^{-1}\lambda(v-f( e^{\us})) 
+4\pi \sum\limits_{j=1}^{m}\delta_{p_j}-4\pi \sum\limits_{k=1}^{n}\delta_{q_k} 
&&\text{ on }M
\\
\label{sys_v}
-\Delta v=&\eps^{-1}\left[\lambda f'(e^{\us})e^{\us}(s-v)
-\eps^{-1}(v-f(e^{\us}))\right]&&\text{ on }M,
\end{align}
where $f:[0,+\infty)\to\R$ is defined by $f(t)=(t-1)/(t+1)$.
In the special case $m=0$ system \eqref{sys_u}--\eqref{sys_v} was introduced in \cite{R2}
with the aim of providing a unified framework for the results in \cite{CK,R1,RT}
and in \cite{CN}. A multiplicity result for \eqref{sys_u}--\eqref{sys_v}
when $m=0$ was obtained in \cite{R3}. 
Our main result concerns the multiplicity of solutions for 
system \eqref{sys_u}--\eqref{sys_v} in the case $m>0$
under the following 
\medskip
\par\noindent
\textbf{Assumptions on $f$:}
\begin{description}
\item[($f0$)] 
$f:[0,+\infty)\to\R$ smooth and $f^{\prime }(t)>0$ $\forall t>0;$
\item[($f1$)]
$f(0)<s<\sup\limits_{t>0}f(t)<+\infty;$
\item[($f2$)] 
$\sup_{t>0}t^4\left\vert f'''(t)\right\vert <+\infty$.
\end{description}
\par
For later use, we note that assumptions $(f0)$--$(f1)$--$(f2)$ imply 
that there exists $f_\infty>s$ such that
\begin{equation}
\label{hypf}
\sup_{t>0}\left[t|f(t)-f_\infty|+t^2f'(t)+t^3|f''(t)|+t^4|f'''(t)|\right]<+\infty.
\end{equation}
Clearly, $f$ defined by $f(t)=(t-1)/(t+1)$ satisfies $(f0)$--$(f1)$--$(f2)$
for every $-1<s<1$.
We restrict our attention to the case $m>n$. It will be clear 
that the case $m<n$ may be treated analogously, while the case $m=n$ requires
an altogether different method and will not be considered here.
Our main result is the following
\begin{Theorem}
\label{thm:main} Let $m>n$ and suppose that $f$ satisfies assumptions 
$(f0)$--$(f1)$--$(f2)$. Then there exists $\lambda_0>0$ with the property that for every fixed 
$\lambda\ge\lambda_0$ there exists $\eps_\lambda>0$ such that for each 
$0<\eps\le\eps_\lambda$ system \eqref{sys_u}--\eqref{sys_v} admits at least \emph{two} solutions.
\end{Theorem}
The remaining part of this note is devoted to the proof of Theorem~\ref{thm:main}.
In Section~\ref{sec:prelims} we prove that system \eqref{sys_u}--\eqref{sys_v} is
equivalent to the following nonlinear elliptic equation of the fourth order:
\begin{align}
\label{fourthintro}
\eps^2&\Delta^{2}u-\Delta u 
\\
=&-\eps\lambda\left[f^{\prime\prime}\left( e^{\sigma +u}\right)
e^{\sigma +u}+f'\left(e^{\sigma +u}\right) \right] e^{\sigma+u}
\left\vert\nabla(\sigma +u)\right\vert ^{2}\nonumber
\\
&+2\eps\lambda\Delta f\left(e^{\sigma+u}\right)
+\lambda^2f'\left(e^{\sigma +u}\right)e^{\sigma +u}
\left(s-f\left(e^{\sigma +u}\right)\right)+4\pi(m-n)\qquad \text{on}\ M,  \notag
\end{align}
where $\sigma$ is the Green function uniquely defined by 
$-\Delta\sigma=4\pi\sum_{j=1}^{m}\delta _{p_j}
-4\pi\sum_{k=1}^{n}\delta _{q_{k}}-4\pi(m-n)$, $\int_M\sigma=0$
(note that $|M|=1$).
By formally setting $\eps=0$ in \eqref{fourthintro} we obtain the 
``limit" equation
\begin{equation}
\label{limiteqintro}
-\Delta u=\lambda^2f'\left(e^{\sigma +u}\right)e^{\sigma +u}
\left(s-f\left(e^{\sigma +u}\right)\right)+4\pi(m-n)
\qquad\text{on }M.
\end{equation}
For $f(t)=(t-1)/(t+1)$ equation \eqref{limiteqintro} describes the vortex solutions for
the $CP(1)$ Chern-Simons model introduced in \cite{KLL1} and analyzed in \cite{CN1}.
When $f(t)=t$ and $s=1$, equation \eqref{limiteqintro} describes the vortex solutions for
the $U(1)$ Chern-Simons model introduced in \cite{HKP,JW}, which has received  considerable attention 
by analysts in recent years, see \cite{CY,DJLW,NT,T}
and the references therein. 
In turn, solutions to \eqref{fourthintro} correspond to 
critical points in the Sobolev space $H^2(M)$ for the functional $\Ie$ 
defined by 
\begin{align*}
\Ie(u)=&\frac{\eps^2}{2}\int\left(\Delta u\right)^2
+\frac{1}{2}\int|\nabla u|^2
+\eps\lambda\int f'\left(e^{\sigma+u}\right)e^{\sigma+u}|\nabla(\sigma+u)|^2
\\
&+\frac{\lambda^2}{2}\int\left(f\left(e^{\sigma +u}\right)-s\right)^2-4\pi(m-n)\int u.
\end{align*}
The two desired solutions for \eqref{sys_u}--\eqref{sys_v}
will be obtained as a local minimum and a mountain pass for $\Ie$
(in the sense of Ambrosetti and Rabinowitz~\cite{AR}). 
The main issue will be to produce a local minimum for $\Ie$. 
To this end, in Section~\ref{sec:locmin} we first construct a supersolution $\up$ for 
equation \eqref{limiteqintro}.
By adapting to the fourth order equation \eqref{fourthintro}
the constrained minimization technique for second order equations
in Brezis and Nirenberg~\cite{BN} (see also Tarantello \cite{T}),
we set $\A=\{u\in H^2(M),\ u\le\up \text{ a.e. on  }M\}$
and we consider $\ue\in\A$ satisfying
$\Ie(\ue)=\min_{\A}\Ie$.
Then $\ue$ is a subsolution for \eqref{fourthintro}.
By an accurate analysis we show that for small values of $\eps$
we have in fact the \textit{strict} inequality $\ue<\up$ everywhere on $M$. 
Consequently, $\ue$ is an internal minimum point for $\Ie$
on $\A$ in the sense of $H^2$, and thus it yields a local minimum for $\Ie$. On the other hand
we have
$\Ie(c)\to-\infty$ on constant functions $c\to+\infty$. Consequently, $\Ie$ has a mountain pass geometry.
In Section~\ref{sec:ps} we prove the Palais-Smale condition for $\Ie$.
At this point, the classical mountain pass theorem in \cite{AR} concludes the proof of
Theorem~\ref{thm:main}. 
The Appendix contains
some simple technical facts which are repeatedly used throughout the proofs.
\par
\medskip\noindent
\textbf{Notation.}
Henceforth, unless otherwise specified, all equations are defined on $M$,
all integrals are taken over $M$ with respect to the Lebesgue measure and
all functional spaces are defined on $M$ in the usual way. In particular, 
we denote by $L^p$, $1\le p\le+\infty$, the Lebesgue spaces and by $H^k$,
$k\ge1$ the Sobolev spaces. We denote by $C>0$ a general constant, 
independent of certain parameters that will be specified in the sequel, 
and whose actual value may vary from line to line.
\section{Preliminaries}
\label{sec:prelims} 
In this section we show that system \eqref{sys_u}--\eqref{sys_v} 
admits a variational structure. We set $A=4\pi (m-n)>0$.
Following a technique introduced by Taubes for self-dual models (see \cite{JT}),
we denote by $\sigma$ the Green's function uniquely defined by
\begin{align*}
&-\Delta\sigma =4\pi \sum\limits_{j=1}^{m}\delta_{p_j}
-4\pi\sum\limits_{k=1}^{n}\delta_{q_k}-A 
\\
&\int\sigma=0
\end{align*}
(recall that $|M|=1$).
Setting $\us=\sigma +u$ system \eqref{sys_u}--\eqref{sys_v} takes
the form
\begin{align}
\label{eq_u_after_taubes} 
-\Delta u&=\eps^{-1}\lambda \left( v-f\left( e^{\sigma+u}\right)\right)+A 
\\
\label{eq_v_after_taubes}
-\Delta v& =\eps ^{-1}\left[ \lambda f^{\prime }\left( e^{\sigma
+u}\right) e^{\sigma +u}(s-v)-\eps ^{-1}\left( v-f\left( e^{\sigma+u}\right)\right)\right]. 
\end{align}
In turn, system \eqref{eq_u_after_taubes}--\eqref{eq_v_after_taubes}
is equivalent to a fourth order equation. We note that equation \eqref{eq_v_after_taubes}
may be written in the equivalent form:
\begin{equation}
\label{vequiv}
-\Delta v+\eps^{-2}(1+\eps\lambda f'(e^{\sigma+u})e^{\sigma+u})v
=\eps^{-2}[f(e^{\sigma+u})+\eps\lambda f'(e^{\sigma+u})e^{\sigma+u}].
\end{equation}
By uniqueness for equation~\eqref{vequiv} for every fixed $u$,
if $v\in L^1$ is a distributional solution for \eqref{vequiv}, then it is in fact $H^1$.
We first show:
\begin{Lemma}
\label{lem:sols are classical}
Suppose $(u,v)\in H^{1}\times H^{1}$ is a weak solution for system 
\eqref{eq_u_after_taubes}--\eqref{eq_v_after_taubes}. Then $(u,v)$ is a classical solution.
\end{Lemma}
\begin{proof}
Throughout this proof, we denote by $\alpha>0$ a general H\"older exponent.
By \eqref{vequiv}, $v\in C^\alpha$. Then by \eqref{eq_u_after_taubes}, $u\in C^{1,\alpha}$.
By Lemma~\ref{lem:lipschitz} and \eqref{hypf}, $f\left(e^{\sigma+u}\right)$ 
and $f'\left(e^{\sigma+u}\right)e^{\sigma+u}$
are Lipschitz continuous. Therefore, by \eqref{vequiv} $v\in C^{2,\alpha}$. 
In turn, by \eqref{eq_u_after_taubes}
$u\in C^{2,\alpha}$ and in particular $(u,v)$ is a classical solution.
\end{proof}
\begin{Lemma}
\label{lem:from_sys_to_fourth}
The couple $(u,v)\in H^1\times H^1$ is a weak
solution for system~\eqref{eq_u_after_taubes}--\eqref{eq_v_after_taubes}
if and only if $u\in H^{2}$ is a weak solution for the fourth order equation
\begin{align}
\label{fourth}
\eps^2\Delta^{2}u-&\Delta u 
\\
=&-\eps\lambda\left[f^{\prime\prime}\left( e^{\sigma +u}\right)
e^{\sigma +u}+f'\left(e^{\sigma +u}\right) \right] 
e^{\sigma+u}\left\vert\nabla(\sigma +u)\right\vert ^{2}\nonumber
\\
&+2\eps\lambda\Delta f\left(e^{\sigma+u}\right)
+\lambda^2f'\left(e^{\sigma +u}\right)e^{\sigma +u}
\left(s-f\left(e^{\sigma +u}\right)\right)+A,  \notag
\end{align}
and $v$ is defined by
\begin{equation}
\label{eq_v_after_taubes_2}
v=-\eps \lambda ^{-1}\Delta u-\eps \lambda ^{-1}A+f\left(
e^{\sigma +u}\right) .  
\end{equation}
\end{Lemma}
\begin{proof} 
Suppose $(u,v)\in H^1\times H^1$ is a weak solution for 
\eqref{eq_u_after_taubes}--\eqref{eq_v_after_taubes}.
Then by Lemma~\ref{lem:sols are classical} we have in particular $u\in H^2$. 
Solving \eqref{eq_u_after_taubes}
for $v$, we obtain \eqref{eq_v_after_taubes_2}.
Inserting the expression for $v$ as given by \eqref{eq_v_after_taubes_2} into
\eqref{eq_v_after_taubes}, we find that $u$ is a distributional solution for
the equation
\begin{align}
\label{step1_equiv_sys_eq} 
\eps^2\Delta^2u-\Delta u  
=&\eps \lambda \Delta f\left( e^{\sigma +u}\right)
+\eps\lambda f'\left(e^{\sigma +u}\right)e^{\sigma +u}\left(\Delta u+A\right)\\
\nonumber
&+\lambda^2f'\left(e^{\sigma+u}\right)e^{\sigma+u}\left(s-f(e^{\sigma +u})\right)+A.  \notag
\end{align}
On the other hand, by the identities \eqref{A} and \eqref{lapl_f} in the Appendix
we have, in the sense of distributions:
\begin{align}
&\Delta f\left( e^{\sigma +u}\right) +f^{\prime }\left( e^{\sigma
+u}\right) e^{\sigma +u}\left(\Delta u+A\right)=\Delta f\left(
e^{\sigma +u}\right) +f^{\prime }\left( e^{\sigma +u}\right) e^{\sigma
+u}\Delta (\sigma +u) \\
&=2\Delta f\left( e^{\sigma +u}\right) -\left\{ f^{\prime \prime }\left(
e^{\sigma +u}\right) e^{\sigma +u}+f^{\prime }\left( e^{\sigma +u}\right)
\right\} e^{\sigma +u}\left\vert \nabla (\sigma +u)\right\vert ^{2}.  \notag
\end{align}
Inserting into \eqref{step1_equiv_sys_eq}, we conclude that $u\in H^2$ satisfies \eqref{fourth}. 
\par
Conversely, suppose $u\in H^2$ is a weak solution for \eqref{fourth}. Then $v$ defined by
\eqref{eq_v_after_taubes_2} belongs to $L^2$, and thus it is a distributional solution for
\eqref{eq_v_after_taubes}. By uniqueness of solutions to 
\eqref{vequiv} for fixed $u$, we conclude that $v\in H^1$. 
\end{proof}
Equation \eqref{fourth} admits a variational formulation, as stated in
the following
\begin{Lemma}
\label{lem:from_fourth_to_functional}$u\in H^2$ is a weak solution for 
\eqref{fourth} if and only if it is a critical point for the $C^1$ functional $\Ie$ 
defined on $H^{2}$ by:
\begin{align}
\Ie(u)=&\frac{\eps^2}{2}\int\left(\Delta u\right)^2
+\frac{1}{2}\int|\nabla u|^2
\label{I_eps(u)} \\
+&\eps\lambda\int f'\left(e^{\sigma+u}\right)e^{\sigma+u}|\nabla(\sigma+u)|^2
+\frac{\lambda^2}{2}\int\left(f\left(e^{\sigma +u}\right)-s\right)^2-A\int u.  
\notag
\end{align}
\end{Lemma}
\begin{proof}
By Lemma~\ref{lem:welldefined} and properties of $f$
as in \eqref{hypf}, $\Ie$ is well-defined and $C^1$ on $H^2$.
We compute, for any $\phi\in H^2$:
\begin{align*}
\left. \frac{d}{dt}\right|_{t=0}\int f'&\left(e^{\sigma+u+t\phi }\right)e^{\sigma+u+t\phi}
\left|\nabla(\sigma+u+t\phi)\right|^2\\
=&\int\left[f''\left(e^{\sigma +u}\right)e^{\sigma+u}
+f'\left(e^{\sigma +u}\right)\right]e^{\sigma +u}\left\vert
\nabla(\sigma+u)\right\vert^2\phi
\\
&+2\int f'\left(e^{\sigma +u}\right)e^{\sigma +u}\nabla(\sigma+u)\cdot\nabla\phi.
\end{align*}
Consequently,
\begin{align}
\label{Iprime}
\bra\Ie'(u),\phi\ket=&\eps^2\int\Delta u\Delta\phi+\int \nabla u\cdot\nabla\phi
\\
\nonumber
&+\eps \lambda \int \left[ f^{\prime \prime}\left(e^{\sigma
+u}\right) e^{\sigma +u}+f^{\prime }\left( e^{\sigma +u}\right)\right]
e^{\sigma +u}\left\vert\nabla(\sigma +u)\right\vert ^{2}\phi  
\\
\nonumber
&+2\eps \lambda \int f^{\prime }\left( e^{\sigma +u}\right)
e^{\sigma +u}\nabla(\sigma +u)\cdot\nabla \phi  
\\
\nonumber
&+\lambda^2\int f'\left(e^{\sigma+u}\right)e^{\sigma+u}
\left(f(e^{\sigma+u})-s\right)\phi-A\int\phi.  
\end{align}
Since
\begin{equation}
\int f'\left(e^{\sigma+u}\right)e^{\sigma +u}\nabla(\sigma+u)\cdot\nabla\phi 
=\int\nabla f\left(e^{\sigma +u}\right)\cdot\nabla \phi 
=-\int\Delta f\left( e^{\sigma +u}\right) \phi,
\end{equation}
it follows that critical points of $\Ie$ correspond to
solutions for \eqref{fourth}, as asserted.
\end{proof}
\section{A local minimum}
\label{sec:locmin} 
Our aim in this section is to prove the existence of a local minimum
for the functional $\Ie$, as stated in the following
\begin{Proposition}
\label{prop:subsol}
There exists $\lambda_0>0$ with the property that for every fixed 
$\lambda\geq\lambda_0$ there exists
$\eps_\lambda>0$, such that for any $0<\eps\le\eps_\lambda$
there exists a solution $\ue$ to \eqref{fourth} corresponding
to a local minimum for the functional $\Ie$.
\end{Proposition}
Throughout this section, we denote by $C>0$ a general constant
independent of $\eps>0$.
Following an idea in \cite{CN1}, we first construct a supersolution for the ``limit" equation
\begin{equation}
\label{limiteq}
-\Delta u=\lambda^2f'\left(e^{\sigma+u}\right)e^{\sigma+u}
\left(s-f\left(e^{\sigma +u}\right)\right)+A,
\end{equation}
which is formally obtained from \eqref{fourth} by setting $\eps=0$.
\begin{Lemma}
\label{lem:exist_supersol}
There exists $\lambda_0>0$ such that
for all $\lambda\ge\lambda_0$ equation \eqref{limiteq}
admits a (distributional) supersolution $\up$.
\end{Lemma}
\begin{proof}
We fix $\rho>0$ such that $B_\rho(p_j)\cap B_\rho(p_l)=\emptyset$
for all $j\neq l$ and
\begin{equation*}
\label{cond_on_rho}
\sum\limits_{j=1}^m|B_\rho(p_j)|<\frac12.  
\end{equation*}
We denote by $g$ a smooth cutoff function satisfying
\begin{equation*}
\label{g+}
g(x)=\begin{cases}
1,\ &\text{if }x\in \bigcup\limits_{j=1}^mB_\rho(p_j)\\
0,&\text{if }x\in M\backslash \bigcup\limits_{j=1}^{m}B_{2\rho}(p_j)
\end{cases}
\end{equation*}
and $0\leq g(x)\leq 1$ for all $x\in M.$
Let $u^{\ast}$ be the function uniquely defined by
\begin{align*}
&-\Delta u^{\ast }=A-4\pi m+8\pi m\left(g-\int g\right)
+4\pi\sum_{k=1}^{n}\delta _{q_k}
\\
&\int u^{\ast }=0.
\end{align*}
We define
\begin{equation*}
\up=u^{\ast}+\bar C,
\end{equation*}
with $\bar C>0$ sufficiently large so that
\begin{equation}
\label{C_big}
f(e^{\sigma +\up})-s>c_0\text{ on }M  
\end{equation}
for some $c_0>0$.
Such a $\bar C$ exists by $(f1)$ since $\sigma+u^\ast$ is bounded below on $M$. 
We claim that for all $\lambda$ sufficiently large,
$\up$ is a supersolution for \eqref{limiteq}.
Indeed, if $x\in \bigcup\limits_{j=1}^mB_\rho(p_j)$, then $g(x)=1$ and
in view of \eqref{C_big}
\begin{align*}
-\Delta \up&\geq A-4\pi m+8\pi m
\left(1-\sum_{j=1}^m\left|B_\rho(p_{j})\right|\right)+4\pi \sum_{j=1}^{n}\delta _{q_k}
\geq A
\\
&\geq \lambda^2f'\left(e^{\sigma+\up}\right)e^{\sigma+\up}
\left(s-f\left(e^{\sigma+\up}\right)\right)+A.
\end{align*}
On the other hand, if $x\in M\setminus \bigcup\limits_{j=1}^mB_\rho(p_j)$, we have:
\begin{align*}
-\Delta \up\geq&A-4\pi m-8\pi m\int g
+4\pi\sum_{k=1}^{n}\delta_{q_k}
\ge A-12\pi m.  
\end{align*}
Let us check that on $M\setminus \bigcup\limits_{j=1}^{m}$ $B_\rho(p_j)$ we have
\begin{equation*}
A-12\pi m\geq \lambda^2f'\left(e^{\sigma+\up}\right)e^{\sigma +\up}
\left(s-f\left(e^{\sigma+\up}\right)\right)+A  
\end{equation*}
for all $\lambda$ sufficiently large.
Indeed, we can choose $C_1>0$ such that
\begin{equation*}
C_1^{-1}\le e^{\sigma +\up}\le C_1\text{ on }M\setminus
\bigcup\limits_{j=1}^{m}B_\rho(p_j).
\end{equation*}
In view of \eqref{C_big} there exists $c_2>0$ such that
\begin{equation*}
f'\left(e^{\sigma+\up}\right)e^{\sigma +\up}\left(s-f\left(e^{\sigma +\up}\right)\right) 
\leq-c_2<0\text{ on }M\setminus\bigcup\limits_{j=1}^{m}B_\rho(p_j).  
\end{equation*}
Therefore, for $\lambda$ large, $\up$ is a subsolution for \eqref{limiteq} in 
$M\setminus\bigcup\limits_{j=1}^{m}B_\rho(p_j)$.
\end{proof}
Henceforth, we fix $\lambda\ge\lambda_0$.
We note that solutions to \eqref{limiteq} correspond to critical points in $H^1$
for the functional $\Io$ defined by
\begin{equation*}
\Io(u)=\frac12\int \left\vert\nabla u\right\vert^{2}
+\frac{\lambda^2}{2}\int(f(e^{\sigma +u})-s)^2-A\int u,
\end{equation*}
for $u\in H^1$.
We define
\begin{equation*}
\A=\left\{ u\in H^{2}:u\leq \up\right\}.
\end{equation*}
Then $\A$ is a convex
closed subset of $H^2$ and consequently there exists $\ue\in H^2$ satisfying
\[
\Ie(\ue)=\min_{\A}\Ie.
\]
Since $\ue-\phi\in\A$ for every $\phi\in H^2$, $\phi\ge0$
we have $\Ie(\ue-\phi)\ge\Ie(\ue)$ for every $\phi\in H^2$, $\phi\ge0$. 
Therefore,
$\ue$ is a weak subsolution for \eqref{fourth}, i.e., it satisfies
\begin{align}
\label{subsol}
\eps^2\Delta^{2}\ue-&\Delta\ue 
\\
\nonumber
\le&-\eps\lambda\left[f^{\prime\prime}\left( e^{\sigma +\ue}\right)
e^{\sigma +\ue}+f'\left(e^{\sigma +\ue}\right) \right] 
e^{\sigma+\ue}\left\vert\nabla(\sigma +\ue)\right\vert ^{2}
\\
\nonumber
&+2\eps\lambda\Delta f\left(e^{\sigma+\ue}\right)
+\lambda^2f'\left(e^{\sigma +\ue}\right)e^{\sigma +\ue}
\left(s-f\left(e^{\sigma +\ue}\right)\right)+A.  
\end{align}
in the weak sense.
The main step towards proving Proposition~\ref{prop:subsol} will be to prove
the \textit{strict} inequality 
$\ue<\up$ on $M$, see Lemma~\ref{lem:u_eps<u_+} below.
We begin by establishing:
\begin{Lemma}
\label{lem:uo}
There exists a subsolution $\uo\in H^1$ for equation \eqref{limiteq} such that
$\ue\weak\uo$ weakly in $H^1$, strongly in $L^p$ for every $p\ge1$ and a.e.\
on $M$. Furthermore, $\uo<\up$.
\end{Lemma}
\begin{proof}
We denote by $\mu\in\A$ the constant function defined by
\[
\mu(x)=\min_M\up\text{ for all }x\in M.
\]
Then 
\begin{equation*}
\label{I_eps(u_eps)<I_eps(m)}
\Ie(\ue)\leq \Ie(\mu)\le C.
\end{equation*}
Since we also have $\int\ue\le\int\up\le C$, we readily derive the estimates
\begin{align}
\label{uefirst}
\eps \left\Vert \Delta \ue\right\Vert _{2}+\left\Vert
\nabla \ue\right\Vert _{2}+\eps \int f^{\prime }(e^{\sigma
+\ue})e^{\sigma +\ue}|\nabla \left( \sigma +u\right)
|^{2}+\left|\int\ue\right|\leq C  
\end{align}
In particular, we have $\|\ue\|_{H^1}\le C$. 
Therefore, by Sobolev embeddings
there exists $\uo\in H^{1}$ such that up to subsequences
$\ue\weak\uo$ weakly in $H^1$, strongly
in $L^p$ for every $p\ge1$ and a.e.\ on $M$.
In particular, $\uo\le\up$ on $M$.
Taking limits in \eqref{subsol}, we find that $\uo$ is a subsolution
for \eqref{limiteq}.
Now the strong maximum principle yields
$\uo<\up$ on $M$.
\end{proof}
Now we can strengthen the convergences stated in Lemma~\ref{lem:uo}.
\begin{Lemma}
\label{lem:eps||_lapl_u_eps||}
The following limits hold:
\begin{itemize}
\item[(i)]
$\lim_{\eps\to0}\Ie(\ue)=\inf_{\A}\Io=\Io(\uo)$
\item[(ii)]
$\lim_{\eps\to0}\eps\left\Vert \Delta \ue\right\Vert_2=0$
\item[(iii)]
$\lim_{\eps\to0}\eps\int f'(e^{\sigma +\ue})e^{\sigma+\ue}|\nabla(\sigma+\ue)|^2=0.$
\end{itemize}
Furthermore, $\uo$ is in fact a solution for \eqref{limiteq}.
\end{Lemma}
\begin{proof}
Proof of (i).
The functional $\Ie$ may be written in the form
\begin{equation}
\label{I_eps_in_terms_of_I_0}
\Ie(u)=\frac{\eps^2}{2}\left\Vert\Delta u\right\Vert_{2}^{2}
+\eps\lambda\int f'(e^{\sigma +u})e^{\sigma+u}
\left\vert\nabla(\sigma +u)\right\vert ^2+\Io(u)
\end{equation}
for every $u\in H^2$.
Consequently
\begin{equation*}
\Ie(\ue)=\inf_{\A}\Ie\geq \inf_{\A}\Io,
\end{equation*}
and therefore
\begin{equation}
\label{liminf_I_eps(u_eps)>inf_I_0}
\liminf\limits_{\eps \rightarrow 0}I_{\eps }(\ue)\geq
\inf_{\A}I_{0}.  
\end{equation}
In order to prove that
\begin{equation}
\label{limsup_I_eps(u_eps)<inf_I_0}
\limsup_{\eps\rightarrow 0}I_{\eps }(\ue)\leq
\inf\limits_{\A}I_{\eps },
\end{equation}
we observe that for any $\eta>0$ we can select $u_{\eta }\in\A$
such that
\begin{equation*}
I_{0}(u_{\eta})\leq\inf_{\A}I_{0}+\eta .
\end{equation*}
Then we have
\begin{align*}
I_{\eps }(\ue)\leq&I_{\eps }(u_{\eta })\leq
I_{0}(u_{\eta })+\circ_\eps(1) 
\\
\leq&\inf_{\A}\Io+\eta +\circ_\eps(1).  
\end{align*}
Therefore
\begin{equation}
\limsup\limits_{\eps\to0}\Ie(\ue)
\leq\inf_{\A}\Io+\eta,
\end{equation}
and since $\eta$ can be chosen arbitrarily small we obtain \eqref{limsup_I_eps(u_eps)<inf_I_0}. 
From \eqref{liminf_I_eps(u_eps)>inf_I_0} and \eqref{limsup_I_eps(u_eps)<inf_I_0} 
we obtain (i).
\par
Proof of (ii)--(iii). Since $\ue\weak\uo$ weakly in $H^1$, we have
\begin{equation*}
\liminf\limits_{\eps\to0}\Io(\ue)\geq\Io(\uo).
\end{equation*}
Therefore,
\begin{align*}
\inf_{\A}\Io=&\lim_{\eps \to0}\Ie(\ue)
\\
=&\lim_{\eps \to0}\left\{\frac{\eps ^{2}}{2}\left\Vert
\Delta \ue\right\Vert _{2}^{2}+\eps \lambda \int f^{\prime
}(e^{\sigma +\ue})e^{\sigma +\ue}\left\vert \nabla
(\sigma +\ue)\right\vert ^{2}+I_{0}(\ue)\right\}  
\\
\geq&I_{0}(u_{0})  
\geq\inf\limits_{\A}I_{0}.  
\end{align*}
Hence,  (ii) and (iii) are established. 
By (i)
we obtain that $u_{\eps}\rightarrow u_{0}$ strongly in 
$H^1$ and $I_{0}(u_{0})=\inf_\A\Io.$
Since we also have $\uo<\up$ (see Lemma~\ref{lem:uo}),
we have that $\uo$ belongs to the interior of $\A$ in the $C^0$-topology.
In particular, $u_{0}$ is a local minimum for $I_{0}$ in the $C^1$-topology. By
the Brezis and Nirenberg argument in \cite{BN}, $\uo$ is a local minimum
for $I_{0}$ in the $H^{1}$-topology and thus it is in fact a solution for 
\eqref{limiteq}.
\end{proof}
Now we are ready to prove the following crucial strict inequality:
\begin{Lemma}
\label{lem:u_eps<u_+}
For every fixed $\lambda\ge\lambda_0$ there exists $\eps_\lambda>0$
such that for every $0<\eps<\eps_\lambda$ there holds
\begin{equation}
\ue<\up\text{ \textsl{on} }M.
\end{equation}
\end{Lemma}
\begin{proof}
We denote
\begin{equation*}
a(u)=-\left[ f^{\prime \prime }(e^{\sigma +u})e^{\sigma
+u}+f^{\prime }(e^{\sigma +u})\right] e^{\sigma +u}\left\vert \nabla \left(
\sigma +u\right) \right\vert ^{2}+2\Delta f(e^{\sigma +u})
\end{equation*}
for all $u\in H^2$ and
\begin{align*}
F_{\eps }=&\eps \lambda a(\ue)+\lambda ^{2}f'(e^{\sigma +\ue})e^{\sigma +\ue}
\left( s-f(e^{\sigma+\ue})\right) +A,  
\\
F_{0}=&\lambda ^{2}f^{\prime }(e^{\sigma +\uo})e^{\sigma+\uo}
\left( s-f(e^{\sigma+\uo})\right)+A.
\end{align*}
Then \eqref{subsol} may be written in the form: 
\begin{equation*}
\eps ^{2}\Delta ^{2}\ue-\Delta \ue\leq F_{\eps }.
\end{equation*}
Now we exploit the decomposition $\eps^2\Delta^2-\Delta=(-\eps^2\Delta+1)(-\Delta)$.
Let $\Ge$ be the Green function for the operator $-\eps^2\Delta+1$ on $M$.
In what follows we shall repeatedly use the properties of $\Ge$ established in 
Lemma~\ref{lem:G_eps} in the Appendix. 
Since $G_{\eps }>0$ on $M$, 
from the above inequality we derive 
\begin{equation}
-\Delta \ue\leq\Ge\ast\Fe.
\label{ineq_for_lapl_u_eps}
\end{equation}
\textbf{Claim}: There exists $1<q<2$ such that
\begin{equation}
\label{|| F_eps_minus_F_0 ||}
\left\Vert F_{\eps}-F_{0}\right\Vert _{q}\rightarrow 0\text{ as }
\eps \rightarrow 0.  
\end{equation}
Proof of \eqref{|| F_eps_minus_F_0 ||}. We only show that 
$\eps\left\Vert a(\ue)\right\Vert _{q}\rightarrow 0$ as $\eps
\rightarrow 0$, since the remaining estimates follow by 
compactness arguments in a straightforward manner.
By identity~\eqref{lapl_f} in the Appendix we may write
\begin{align*}
a(\ue)=&\left[f''(e^{\sigma+\ue})e^{\sigma +u_{\eps}}
+f^{\prime}(e^{\sigma +u_{\eps}})\right]
e^{\sigma+u_{\eps}}\left\vert\nabla(\sigma+u_{\eps})\right\vert^2 
\\
&+2f'(e^{\sigma+\ue})e^{\sigma +\ue}\Delta
(\sigma +\ue).  
\end{align*}
Therefore, in view of \eqref{hypf} it suffices to show that as $\eps\rightarrow 0$
\begin{align}
\label{a_eps_2.1}
&\eps\left\Vert\left[f''(e^{\sigma+u_{\eps}})e^{\sigma +\ue}+f'(e^{\sigma +\ue})\right]
e^{\sigma+\ue}\left|\nabla \sigma \right|^2\right\Vert _{q}\rightarrow 0
\\ 
\label{epsgrad}
&\eps\|\,|\nabla\ue|^2\|_q\to0
\\
\label{a_eps_2.2}
&\eps \left\Vert f'(e^{\sigma +\ue})e^{\sigma+\ue}
\Delta (\sigma+\ue)\right\Vert _{q}\to0.  
\end{align}
To see \eqref{a_eps_2.1}, note that by Lemma~\ref{lem:welldefined} in the Appendix 
and properties of $f$ as stated in \eqref{hypf},
\begin{align*}
\left[f''(e^{\sigma+\ue})e^{\sigma+\ue}+f'(e^{\sigma +\ue})\right]
e^{\sigma+\ue}\left|\nabla\sigma \right|^2\le C\left(1+e^{\ue}+e^{-\ue}\right).
\end{align*}
By \eqref{uefirst} and the Moser-Trudinger inequality (see, e.g., Aubin~\cite{Aubin}), 
we have
$\|e^{\ue}\|_q\le Ce^{\int\ue}e^{q\|\nabla\ue\|_2^2}\le C$.
Similarly, we obtain $\|e^{-\ue}\|_q\le C$. Therefore, \eqref{a_eps_2.1}
is established.
To see \eqref{epsgrad}, let $1<q<\alpha<2$.
Then by H\"older's inequality and \eqref{uefirst},
\begin{align*}
\int|\nabla\ue|^{2q}=&\int|\nabla\ue|^\alpha|\nabla\ue|^{2q-\alpha}
\\
\le&\left(\int|\nabla\ue|^2\right)^{\alpha/2}
\left(\int|\nabla\ue|^{2(2q-\alpha)/(2-\alpha)}\right)^{(2-\alpha)/2}
\\
\le&C\|\nabla\ue\|_{2(2q-\alpha)/(2-\alpha)}^{2q-\alpha}.
\end{align*}
Consequently, in view of Lemma~\ref{lem:eps||_lapl_u_eps||}--(ii):
\begin{align*}
\eps\|\,|\nabla\ue|^2\|_q\le&C\eps\|\nabla\ue\|_{2(2q-\alpha)/(2-\alpha)}^{(2q-\alpha)/q}
\\
\le&C\eps^{1-(2q-\alpha)/q}\left(\eps\|\nabla\ue\|_{2(2q-\alpha)/(2-\alpha)}
\right)^{(2q-\alpha)/q}
\\
\le&C\eps^{1-(2q-\alpha)/q}\|\eps\Delta\ue\|_2^{(2q-\alpha)/q}\to0,
\end{align*}
and \eqref{epsgrad} follows.
Finally, by identity~\eqref{A} in the Appendix and \eqref{hypf} we have
\begin{align*}
\eps\Vert f'(e^{\sigma+\ue})e^{\sigma+\ue}&\Delta(\sigma+\ue)\Vert_{q}
\\
=&\eps\left\Vert Af^{\prime }(e^{\sigma +\ue})e^{\sigma+\ue}
+f'(e^{\sigma +\ue})e^{\sigma +\ue}\Delta\ue\right\Vert_{q}
\\
\leq&\eps \left\Vert Af'(e^{\sigma+\ue})e^{\sigma+\ue}\right\Vert_{q}
+\eps\left\Vert f'(e^{\sigma+\ue})e^{\sigma+\ue}\Delta\ue\right\Vert_{q}
\\
=&\eps\left\Vert f^{\prime }(e^{\sigma+\ue})e^{\sigma +\ue}\Delta\ue\right\Vert_q
+\circ_\eps(1).
\end{align*}
By the H\"older inequality, Lemma~\ref{lem:eps||_lapl_u_eps||}-(ii) and \eqref{hypf} we have
\begin{gather*}
\eps \left\Vert f^{\prime }(e^{\sigma +\ue})e^{\sigma
+\ue}\Delta \ue\right\Vert _{q}
\leq C\eps\left\Vert\Delta\ue\right\Vert_{2}\to0\text{ as }\eps\to0, 
\end{gather*}
which yields \eqref{a_eps_2.2}. We conclude that $\eps\|a(\ue)\|_q\to0$,
as asserted, and the desired claim \eqref{|| F_eps_minus_F_0 ||}
follows.
\par
By \eqref{|| F_eps_minus_F_0 ||}
and properties of $\Ge$ as in Lemma~\ref{lem:G_eps} we derive
\begin{align}
\label{G_eps_star_F_eps_minus_F_0} 
\left\Vert\Ge\ast F_{\eps }-F_{0}\right\Vert_q
\leq&\left\Vert\Ge\ast\left( F_{\eps }-F_{0}\right) \right\Vert_q
+\left\Vert\Ge\ast F_{0}-F_{0}\right\Vert_q
\\
\nonumber
\leq&\left\Vert F_{\eps }-F_{0}\right\Vert_q
+\eps^2\left\Vert\Delta F_{0}\right\Vert_q\rightarrow 0  
\end{align}
as $\eps \rightarrow 0.$ We define $w_{\eps }$ as the unique
solution for
\begin{equation*}
(-\Delta +1)w_{\eps }=\Ge\ast F_{\eps }+u_{_{\eps}}.
\end{equation*}
Then in view of \eqref{ineq_for_lapl_u_eps} we have
\begin{equation*}
(-\Delta +1)\left( \ue-w_{\eps }\right) \leq 0
\end{equation*}
and therefore by the maximum principle
\begin{equation*}
\ue\leq w_{\eps },\text{ on }M.
\end{equation*}
Since $u_{0}$ satisfies \eqref{limiteq}, we have
\begin{equation*}
-\Delta u_{0}=F_{0}.
\end{equation*}
Consequently,
\begin{equation*}
(-\Delta +1)\left( w_{\eps }-u_{0}\right)=G_{\eps }\ast
F_{\eps }-F_{0}+\ue-u_{0},
\end{equation*}
and therefore \eqref{G_eps_star_F_eps_minus_F_0}, Lemma~\ref{lem:uo}
and standard elliptic estimates yield
\begin{equation*}
\left\Vert w_{\eps }-u_{0}\right\Vert _{C^\alpha}
\leq C\left(\left\Vert G_{\eps }\ast F_{\eps }-F_{0}\right\Vert _{q}+\left\Vert
\ue-u_{0}\right\Vert _{q}\right)\to0.  
\end{equation*}
In particular, $w_{\eps }$ converges uniformly to $u_{0}$. 
Taking into account that $\ue\leq w_{\eps }$ and 
$u_{0}<\up$ on $M$, we conclude that for all $\eps>0$ 
sufficiently small we have the desired strict inequality $\ue<\up$.
\end{proof}
Now we can provide the
\begin{proof}[Proof of Proposition~\ref{prop:subsol}] 
Let $\lambda_0>0$ as in Lemma~\ref{lem:exist_supersol} and for every fixed $\lambda\ge\lambda_0$
let $\eps_\lambda>0$ as in Lemma~\ref{lem:u_eps<u_+}. Then by Lemma~\ref{lem:u_eps<u_+}
the function $\ue$ defined by $\Ie(\ue)=\min_\A\Ie$ satisfies the strict inequality
$\ue<\up$ for every
$0<\eps<\eps_\lambda$.
In particular, by the Sobolev embedding 
$\left\Vert u\right\Vert_{\infty }\leq C\left\Vert u\right\Vert_{H^2}$
for all $u\in H^2$,
for every $\eps>0$ sufficiently small there exists an
$H^2$-neighborhood of $\ue$ entirely
contained in $\A$. Therefore, for such values of $\eps$, 
$\ue$ belongs to the interior of $\A$ in the sense of $H^2$.
It follows that $\ue$ is a critical point for $\Ie$ corresponding to a local minimum, 
as asserted.
\end{proof}
\section{The Palais-Smale condition}
\label{sec:ps}
In this section we prove the Palais-Smale condition for $\Ie$
for every fixed $\eps,\lambda>0$. 
\begin{Proposition}
\label{prop:ps}
For every fixed $\eps,\lambda>0$ the functional $I_{\eps }$
satisfies the Palais-Smale condition.
\end{Proposition}
We denote by $(\uj)$, $\uj\in H^2$, $j=1,2,3,\dots$ a 
Palais-Smale sequence for the functional
$I_{\eps }$. That is, $(\uj)$ satisfies:
\begin{align}
\label{P_S_1}
&\Ie(u_{j})\rightarrow \alpha \in \mathbb{R}, 
\\
\label{P_S_2}
&\left\Vert\Ie'(u_{j})\right\Vert _{H^{-1}}\rightarrow 0
\end{align}
as $j\to+\infty$.
We have to show that $(\uj)$ admits a subsequence strongly convergent in $H^2$.
By compactness, it suffices to show that $(u_{j})$ is bounded in $H^{2}.$ It
will be useful to decompose $u_{j}$ in the following way
\begin{equation*}
u_{j}=u_{j}^{\prime }+c_{j},\text{ where }\int u_{j}^{\prime }=0\text{ and }
c_{j}\in \mathbb{R}.  
\end{equation*}
Then condition \eqref{P_S_1} is equivalent to
\begin{align}
\label{I(uj) ->alpha} 
\Ie(u_{j})=&\frac{\eps ^{2}}{2}\int (\Delta u_{j})^{2}+
\frac12\int\left\vert\nabla u_{j}\right\vert ^{2}
+\eps\lambda\int f^{\prime }(e^{\sigma +u_{j}})e^{\sigma+u_{j}}
\left\vert\nabla(\sigma+u_{j})\right\vert^{2}
\\
\nonumber
&+\frac{\lambda ^{2}}{2}\int (f(e^{\sigma +u_{j}})-s)^{2}
-Ac_{j}\rightarrow \alpha  
\end{align}
and (\ref{P_S_2}) implies (see \eqref{Iprime})
\begin{align}
\label{<i',uj'>} 
\circ_{j}(1)&\left\Vert \Delta u_{j}\right\Vert _{2}
=\bra I_{\eps}^{\prime }(u_{j}),u_{j}^{\prime }\ket
=\eps^{2}\int(\Delta u_{j})^{2}+\int \left\vert \nabla u_{j}\right\vert ^{2}
\\
\nonumber
&+\eps\lambda \int 
\left[ f^{\prime \prime }(e^{\sigma +u_{j}})e^{\sigma+u_{j}}
+f^{\prime }(e^{\sigma +u_{j}})\right] e^{\sigma +u_{j}}\left\vert
\nabla (\sigma +u_{j})\right\vert ^{2}u_{j}^{\prime }
\\
\nonumber
&+2\eps\lambda\int f'(e^{\sigma+\uj})e^{\sigma+\uj}\nabla(\sigma+\uj)\cdot\nabla\uj
\\
\nonumber
&+\lambda^{2}\int f^{\prime }(e^{\sigma +u_{j}})e^{\sigma+u_{j}}(f(e^{\sigma +u_{j}})-s)u_{j}'
\end{align}
It is readily checked that $\cj\ge-C$ for some $C>0$. Indeed, by \eqref{I(uj) ->alpha} 
we have
\begin{equation}
\label{||_lapl_uj_||-cj}
-Ac_{j}\leq \frac{\eps ^{2}}{2}\int (\Delta
u_{j})^{2}-Ac_{j}\leq I_{\eps }(u_{j})\leq C.
\end{equation}
Furthermore, if either $c_j\leq C$ or $\left\Vert \Delta u_{j}\right\Vert _{2}\leq
C, $ then $u_{j}$ is bounded in $H^{2}$. 
Indeed, if $c_{j}\leq C$ then we readily obtain from \eqref{||_lapl_uj_||-cj}
that $\left\Vert \Delta u_{j}\right\Vert _{2}$
$\leq C$.
Suppose $\left\Vert \Delta u_{j}\right\Vert_{2}\leq C.$ 
Then by Sobolev embeddings we also have 
$\int\left\vert\nabla u_{j}\right\vert ^{2}\leq C$ and 
$\left\Vert u_{j}^{\prime }\right\Vert_\infty\leq C.$ We have
\begin{align*}
\int f'&(e^{\sigma+\uj})e^{\sigma+\uj}|\nabla(\sigma+\uj)|^2 
\\
=&\int\nabla(f(e^{\sigma+\uj}))\cdot\nabla(\sigma+\uj)
=-\int f(e^{\sigma+\uj})\Delta(\sigma+\uj)
\\
=&-A\int f(e^{\sigma+\uj})+4\pi mf_\infty-4\pi nf(0)-\int f(e^{\sigma+\uj})\Delta\uj  
\end{align*}
and therefore by \eqref{hypf}
\[
\int f'(e^{\sigma+\uj})e^{\sigma+\uj}|\nabla(\sigma+\uj)|^2
\le C(1+\|\Delta\uj\|_2)\le C.
\]
On the other hand the term $\int (f(e^{\sigma+\uj})-s)^2$ is bounded. 
Therefore we derive from \eqref{I(uj) ->alpha} that
\[
\alpha+\circ_j(1)=\Ie(\uj)\le-A\cj+C
\]
and consequently $\cj\le C$.
In view of the above remarks, henceforth we assume that
\begin{equation}
\label{uj_and_cj_to_infinity}
\left\Vert \Delta u_{j}\right\Vert _{2}\to+\infty 
\text{ and }
c_{j}\to+\infty
\text{ as }j\to+\infty.
\end{equation}
By \eqref{||_lapl_uj_||-cj} and assumption
\eqref{uj_and_cj_to_infinity} we then have
\begin{equation}
\label{Deltauc}
\left\Vert \Delta u_{j}\right\Vert _{2}\leq Cc_{j}^{1/2}.
\end{equation}
The following identity will be useful.
\begin{Lemma}
\label{lem:int_by_parts}
For all $u\in H^2$ the following identity holds:
\begin{align*}
\int&\left[ f''(e^{\sigma+u})e^{\sigma+u}
+f'(e^{\sigma+u})\right]e^{\sigma+u}\left\vert\nabla\left(\sigma+u\right)\right\vert^2u   
\\
&\quad+2\int f'(e^{\sigma+u})e^{\sigma+u}\nabla\left(\sigma+u\right)\cdot\nabla u  
\\
=&\int f'(e^{\sigma+u})e^{\sigma+u}\nabla\left(\sigma+u\right)\cdot\nabla u
-\int f'(e^{\sigma+u})e^{\sigma +u}\Delta\left(\sigma+u\right)u.  
\end{align*}
\end{Lemma}
\begin{proof}
Integrating by parts we have
\begin{align}
\int [f''&(e^{\sigma+u})e^{\sigma +u}+f'(e^{\sigma+u})]e^{\sigma+u}
\left\vert\nabla\left(\sigma+u\right)\right\vert^2u 
\\
\nonumber
=&\int \nabla \left[ f^{\prime }(e^{\sigma +u})e^{\sigma +u}\right]\cdot\nabla
\left(\sigma+u\right) u  
\\
\nonumber
=&-\int f'(e^{\sigma +u})e^{\sigma+u}\Delta\left(\sigma+u\right)u
-\int f'(e^{\sigma +u})e^{\sigma +u}\nabla \left(\sigma+u\right)\cdot\nabla u.  
\end{align}
The asserted identity follows.
\end{proof}
Now we can provide the 
\begin{proof}[Proof of Proposition~\ref{prop:ps}]
By \eqref{<i',uj'>} and Lemma~\ref{lem:int_by_parts} we have
\begin{align*}
\circ_j(1)&\|\Delta u_j\|_2
\\
\nonumber
\geq&\eps ^{2}\int(\Delta u_{j})^{2}
+\eps \lambda \int[ f^{\prime \prime }(e^{\sigma +u_{j}})e^{\sigma
+u_{j}}
+f^{\prime }(e^{\sigma +u_{j}})] e^{\sigma +u_{j}}\left\vert
\nabla (\sigma+u_{j})\right\vert ^{2}u_{j}^{\prime }  
\\
\nonumber
+&2\eps \lambda \int f'(e^{\sigma +u_j})e^{\sigma+u_j}\nabla(\sigma+u_j)\cdot\nabla u_j
\\
\nonumber
&+\lambda^2\int f'(e^{\sigma +u_j})e^{\sigma +u_j}(f(e^{\sigma +u_j})-s)\ujp
\\
\nonumber
=&\eps ^{2}\left\Vert \Delta u_{j}\right\Vert _{2}^{2}
+\eps\lambda\int f^{\prime }(e^{\sigma +u_{j}})e^{\sigma +u_{j}}\nabla
\left( \sigma +u_{j}\right)\cdot\nabla u_{j}   
\\
\nonumber
&-\eps\lambda\int f^{\prime }(e^{\sigma +u_{j}})e^{\sigma +u_{j}}\Delta \left(
\sigma +u_{j}\right) u_{j}^{\prime }
\\
\nonumber
&+\lambda^{2}\int f^{\prime
}(e^{\sigma +u_{j}})e^{\sigma +u_{j}}(f(e^{\sigma +u_{j}})-s)u_{j}^{\prime}
\end{align*}
and therefore, since 
$\int f'(e^{\sigma +u_j})e^{\sigma+u_j}|\nabla u_j|^2\ge0$,
\begin{align}
\label{littleo}
\circ_j(1)&\|\Delta u_j\|_2\ge\eps^2\|\Delta\uj\|_2^2
+\eps\lambda\int f^{\prime }(e^{\sigma +u_{j}})e^{\sigma +u_{j}}
\nabla\sigma\cdot\nabla u_{j} 
\\
\nonumber  
-&\eps\lambda\int f^{\prime }(e^{\sigma +u_{j}})e^{\sigma +u_{j}}\Delta \left(
\sigma +u_{j}\right) u_{j}^{\prime }
\\
\nonumber
&+\lambda^2\int f'(e^{\sigma+\uj})e^{\sigma+\uj}\left(f(e^{\sigma+\uj})-s\right)\ujp.
\end{align}
By properties of $f$  and the Sobolev embedding 
$\left\Vert u_{j}^{\prime }\right\Vert _{\infty }\leq C\left\Vert \Delta
u_{j}\right\Vert _{2}$ we have
\begin{equation}
\label{f'(t)t(f(t)-s)u'j}
\left\vert \int f^{\prime }(e^{\sigma +u_{j}})e^{\sigma +u_{j}}(f(e^{\sigma
+u_{j}})-s)u_{j}^{\prime }\right\vert \leq C\left\Vert \Delta
u_{j}\right\Vert _{2}.  
\end{equation}
By the H\"older inequality and Sobolev embeddings we have
\begin{align}
\label{f'(t)tDsigmaDuj}
\left\vert \int f^{\prime }(e^{\sigma +u_{j}})e^{\sigma +u_{j}}
\nabla\sigma\cdot\nabla u_{j}\right\vert
\le C\|\nabla\sigma\|_p\|\nabla\uj\|_{p'}
\le C\left\Vert \Delta u_{j}\right\Vert _{2},  
\end{align}
for any $1\leq p<2.$
By \eqref{A} in the Appendix and Sobolev embeddings we have
\begin{align}
\label{Apalais}
\left\vert \int f^{\prime }(e^{\sigma +u_{j}})e^{\sigma +u_{j}}\right.
\Delta\sigma\,u_{j}^{\prime }\left. \right\vert =A\left\vert \int
f^{\prime }(e^{\sigma +u_{j}})e^{\sigma +u_{j}}u_{j}^{\prime }\right\vert
\leq C\left\Vert \Delta u_{j}\right\Vert_{2}.
\end{align}
Finally, we claim that
there exists $\bar\jmath$ such
that for all $j\ge\bar\jmath$
\begin{equation}
\label{epspalais}
\left\vert\int f'(e^{\sigma +u_{j}})e^{\sigma +u_{j}}\Delta
u_{j}u_{j}\right\vert \leq \frac{\eps^2}{2}\left\Vert \Delta
u_{j}\right\Vert_2^2.
\end{equation}
To prove \eqref{epspalais}, we write for $\rho>0$  
\begin{align*}
\int f^{\prime }(e^{\sigma +u_{j}})e^{\sigma +u_{j}}\Delta u_{j}u_{j} 
=&\int_{\cup _{h=1}^{n}B_\rho(q_{j})}f^{\prime }(e^{\sigma
+u_{j}})e^{\sigma +u_{j}}\Delta u_{j}u_{j}^{\prime } 
\\
&+\int_{M\backslash \left[ \cup _{h=1}^{n}B_\rho(q_{j})\right]}
f^{\prime}(e^{\sigma +u_{j}})e^{\sigma +u_{j}}\Delta u_{j}
u_{j}^{\prime }.  
\end{align*}
In view of the assumptions on $f$, the Cauchy-Schwarz inequality and Sobolev embeddings, 
we estimate
\begin{align*}
\Big|\int_{\cup_{h=1}^nB_\rho(q_j)}&f'(e^{\sigma+\uj})e^{\sigma+\uj}\Delta u_{j}\ujp\Big|
\\
\leq&C\left(\int_{\cup_{h=1}^{n}B_\rho(q_{j})}(\Delta u_{j})^{2}\right)^{1/2}
\left(\int_{\cup _{h=1}^{n}B_\rho(q_{j})}\ujp^2\right)^{1/2}
\\
\leq&C\left\Vert\Delta u_{j}\right\Vert_2\left\Vert\uj^{\prime }\right\Vert_\infty
\left(\sum_{h=1}^{n}\left\vert B_\rho(q_{j})\right\vert\right)^{1/2}
\leq C\rho \left\Vert \Delta u_{j}\right\Vert_{2}^{2}.
\end{align*}
Therefore, we may choose $\rho>0$ such that
\begin{equation}
\label{rho}
\left|\int_{\cup _{h=1}^{n}B_\rho(q_{j})}f^{\prime }(e^{\sigma +u_{j}})e^{\sigma
+u_{j}}\Delta u_{j}u_{j}^{\prime}\right|
\le\frac{\eps^2}{4}\|\Delta\uj\|_2^2.
\end{equation}
We define
\begin{equation*}
e_0=\min_{M\backslash \left[ \cup _{h=1}^{n}B_\rho(q_{j})\right]}e^{\sigma }>0.
\end{equation*}
By \eqref{uj_and_cj_to_infinity}, \eqref{Deltauc} and the embedding $\|\ujp\|_\infty\le C\|\Delta\uj\|_2$
we have
\begin{align*}
\min_{M\backslash \left[ \cup _{h=1}^{n}B_\rho(q_{j})\right]}e^{\sigma +u_{j}}
\ge&\min_{M\backslash \left[ \cup _{h=1}^{n}B_\rho(q_{j})\right]}
e^{\sigma-\|u_{j}^{\prime }\|_\infty+c_{j}}
\geq e_0e^{-C\sqrt{c_{j}}+c_{j}}
\\
\ge&e_0e^{\cj/2}\to+\infty
\text{ as }j\to+\infty.
\end{align*}
Therefore, by properties of $f$ and since $c_{j}\rightarrow \infty$, for every $\mu>0$
there exists $j_\mu\in\mathbb{N}$ such that if $j\ge j_\mu$ then
$f^{\prime }(e^{\sigma+u_{j}})e^{\sigma +u_{j}}\leq\mu $
on $M\setminus\left[\cup _{h=1}^{n}B_\rho(q_{j})\right]$.
We conclude that for $j\ge j_\mu$ we have
\begin{equation*}
\left\vert \int_{M\backslash \left[ \cup _{h=1}^{n}B_\rho(q_{j})\right]
}f^{\prime }(e^{\sigma +u_{j}})e^{\sigma +u_{j}}\Delta u_{j}
u_{j}^{\prime }\right\vert \leq C\mu\left\Vert \Delta
u_{j}\right\Vert _{2}^{2}.
\end{equation*}
We choose $\mu>0$ such that
\begin{equation}
\label{mu}
\left\vert \int_{M\backslash \left[ \cup _{h=1}^{n}B_\rho(q_{j})\right]
}f^{\prime }(e^{\sigma +u_{j}})e^{\sigma +u_{j}}\Delta u_{j}
u_{j}^{\prime }\right\vert \leq \frac{\eps^2}{4}\|\Delta\uj\|_2^2.
\end{equation}
Now \eqref{rho} and \eqref{mu} yield \eqref{epspalais} with $\bar\jmath=j_\mu$.
\par
Now we can conclude the proof of Proposition \ref{prop:ps}. 
Indeed, inserting the estimates
\eqref{f'(t)t(f(t)-s)u'j}--\eqref{f'(t)tDsigmaDuj}--\eqref{Apalais}--\eqref{epspalais}
into \eqref{littleo} we obtain 
\begin{equation*}
\eps ^{2}\left\Vert \Delta u_{j}\right\Vert _{2}^{2}\leq
C\left\Vert \Delta
u_{j}\right\Vert _{2}+\frac{\eps^2}{2}\left\Vert \Delta u_{j}\right\Vert_{2}^{2}
\end{equation*}
and consequently we derive that $\|\Delta\uj\|_2\le C$.
This is a contradiction since we have assumed \eqref{uj_and_cj_to_infinity}.
\end{proof}
Now we can finally prove our main result:
\begin{proof}[Proof of Theorem~\ref{thm:main}]
By Proposition \ref{prop:subsol}, 
there exists $\lambda_0>0$ such that for every $\lambda\ge\lambda_0$
fixed, there exists $\eps_\lambda>0$ such that for every
$0<\eps\le\eps_\lambda$
the functional $\Ie$ admits a
critical point corresponding to a local minimum. By
Proposition~\ref{prop:ps}, $\Ie$ satisfies the Palais-Smale
condition. If $\ue$ is not a strict local minimum, it is known
that $\Ie$ has a continuum of critical points (see, e.g.,
\cite{T}). In particular, $\Ie$ has at least two critical points.
If $\ue$ is a strict local minimum, we note that on constant functions
$c\to+\infty$ we have $\Ie(c)\to-\infty$. Therefore $\Ie$ admits
a mountain pass structure in the sense of Ambrosetti and
Rabinowitz \cite{AR}. Hence by the mountain pass theorem
\cite{AR} we obtain the existence of a second critical point for
$\Ie$. In either case, we conclude that the fourth order equation
\eqref{fourth} admits at least two solutions.
By the equivalences as stated in Lemma~\ref{lem:sols are classical} and in 
Lemma~\ref{lem:from_sys_to_fourth},
system \eqref{eq_u_after_taubes}--\eqref{eq_v_after_taubes} 
admits at least two solutions,
as asserted.
\end{proof}
\section{Appendix}
\label{sec:appendix}
We collect in this Appendix the proofs of some simple properties which have been 
repeatedly used 
throughout this note. Recall that $\sigma$ is defined as the unique 
distributional solution for $-\Delta\sigma=4\pi\sum_{j=1}^m\delta_{p_j}
-4\pi\sum_{k=1}^m\delta_{q_k}$, $\int\sigma=0$.
Therefore, there exist smooth functions $\gamma_j$, $\theta_k$
and $\rho>0$
such that $\sigma(x)=\gamma _j(x)+\log \vert x-p_{j}\vert^{-2}$
in $B_\rho(p_{j})$ and $\sigma(x)=\theta_k(x)+\log\vert x-q_k\vert^2$
in $B_\rho(q_{k})$.
\begin{Lemma}
\label{lem:welldefined}
Suppose $\phi:[0,+\infty)\to\R$ is a smooth function satisfying
\[
|\phi(t)|\le C_\phi\min\{t,t^{-1}\}
\]
for some $C_\phi>0$.
Then there exists $\overline C_\phi>0$ depending  on $\phi$ only
such that for all measurable functions $u$ we have
\[
\phi(e^{\sigma+u})|\nabla\sigma|^2\le\overline C_\phi(1+e^u+e^{-u}).
\]
\end{Lemma} 
\begin{proof}
In $M\setminus\cup_{j=1}^mB_\rho(p_{j})\setminus\cup_{k=1}^nB_\rho(q_{k})$
we have
\[
\phi(e^{\sigma+u})|\nabla\sigma|^2\le C_\phi
\sup_{M\setminus\cup_{j=1}^mB_\rho(p_{j})\setminus\cup_{k=1}^nB_\rho(q_{k})}
|\nabla\sigma|^2.
\]
In $B_\rho(p_j)$ we have
\begin{align*}
\phi(e^{\sigma+u})|\nabla\sigma|^2\le C_\phi e^{-(\sigma+u)}|\nabla\sigma|^2
\le C_\phi\sup_{B_\rho(p_j)}(e^{-\sigma}|\nabla\sigma|^2)\,e^{-u}.
\end{align*}
In $B_\rho(q_k)$ we have
\begin{align*}
\phi(e^{\sigma+u})|\nabla\sigma|^2\le C_\phi e^{\sigma+u}|\nabla\sigma|^2
\le C_\phi\sup_{B_\rho(q_k)}(e^\sigma|\nabla\sigma|^2)\,e^u.
\end{align*}
Now the asserted estimate follows.
\end{proof}
\begin{Lemma}
\label{lem:lipschitz}
If $u\in C^1$ and $\phi:[0,+\infty)\to\R$ is a smooth function
satisfying 
\[
|\phi'(t)|\le C_\phi\min\{t^{-1/2},t^{-3/2}\}
\]
for some $C_\phi>0$, 
then $\phi(e^{\sigma+u})$ is Lipschitz continuous on $M$
(with Lipschitz constant depending on $u$).
\end{Lemma}
\begin{proof}
We need only check the claim near the vortex points $p_j$, $q_k$. 
By the mean value theorem we have, for $x,y\in B_\rho(p_j)$, $x,y\neq p_j$:
\begin{align*}
\phi&(e^{(\sigma+u)(x)})-\phi(e^{(\sigma+u)(y)})\\
=&\phi'(e^{(\sigma+u)(y+\theta(x-y))})e^{(\sigma+u)(y+\theta(x-y))}\nabla(\sigma+u)(y+\theta(x-y))
\cdot(x-y)
\end{align*}
for some $0\le\theta\le1$.
By properties of $\sigma$ and $\phi$, 
\begin{align*}
|\phi&(e^{(\sigma+u)(x)})-\phi(e^{(\sigma+u)(y)})|
\\
\le&C_\phi[e^{-(\sigma+u)(y+\theta(x-y))}]^{1/2}|\nabla(\sigma+u)(y+\theta(x-y))|\,|x-y|
\\
\le&C_\phi e^{\|u\|_\infty/2}|y+\theta(x-y)-p_j|
\left(1+|y+\theta(x-y)-p_j|^{-1}+\|\nabla u\|_\infty\right)
\,|x-y|
\\
\le&C_\phi e^{\|u\|_\infty/2}(1+\|\nabla u\|_\infty)\,|x-y|.
\end{align*}
A similar argument yields Lipschitz continuity near the $q_k$'s, and the statement follows.
\end{proof}
\begin{Lemma}
\label{lem:id}
For any $u\in H^2$ the following identities hold, in the sense of distributions:
\begin{equation}
\label{A}
f'\left(e^{\sigma+u}\right)e^{\sigma +u}\Delta\sigma 
=Af'\left(e^{\sigma+u}\right)e^{\sigma +u}  
\end{equation}
and
\begin{align}
\label{lapl_f}
\Delta&f\left(e^{\sigma+u}\right)\\
\nonumber
=&\left\{f''\left(e^{\sigma+u}\right)e^{\sigma+u}+f'\left(e^{\sigma+u}\right)\right\}e^{\sigma+u}
|\nabla(\sigma+u)|^2
+f'\left(e^{\sigma+u}\right)e^{\sigma+u}\Delta(\sigma+u).  
\end{align}
\end{Lemma}
\begin{proof}
By \eqref{hypf} the function $f'\left(e^{\sigma+u}\right)e^{\sigma +u}$ may be extended
by continuity to the whole of $M$ by setting it equal to $0$ at $p_j$, $q_k$. In view 
of the definition of $\sigma$ we obtain \eqref{A}.
\par
Since $u\in H^2$, \eqref{lapl_f} holds pointwise almost everywhere on $M$.
By Lemma~\ref{lem:welldefined}, the right hand side of \eqref{lapl_f} belongs to $L^1$,
and it is absolutely continuous in the $x$-variable, for almost every fixed $y$.
Therefore, \eqref{lapl_f} holds in the sense of distributions.
\end{proof}
Finally, we prove some properties for the Green function $\Ge$ for
the operator $-\eps^2\Delta+1$ on $M$.
\begin{Lemma}
\label{lem:G_eps} 
Let $\Ge=\Ge(x,y)$ be the
Green function defined by
\begin{equation*}
(-\eps^2\Delta _x+1)\Ge=\delta_y\ \text{on}\ M.
\end{equation*}
Then
\begin{itemize}
\item[i)] 
$G_{\eps }>0$ on $M\times M$ and for every fixed
$y\in M$ we have $G_{\eps }\rightharpoondown \delta _{y}$ as
$\eps \rightarrow 0, $ weakly in the sense of measures;
\item[ii)] 
$\left\Vert G_{\eps }\ast h\right\Vert _{q}\leq
\left\Vert h\right\Vert _{q}$ for all $1\leq q\leq+\infty$;
\item[iii)] 
If $\Delta h\in L^q$ for some $q\geq 1$ then
$\left\Vert G_{\eps }\ast h-h\right\Vert_q\leq
\eps ^{2}\left\Vert \Delta h\right\Vert_q$.
\end{itemize}
\end{Lemma}
\begin{proof}
Proof of (i). Note that since $-\eps^2\Delta+1$ is coercive, $\Ge$
is well defined (e.g., by Stampacchia's duality argument
\cite{St}). By the maximum principle, $\Ge>0$ on $M\times M$.
Integrating over $M$ with respect to $x$, we have
$\int\Ge(x,y)\dx=\int|\Ge(x,y)|\dx=1$ and therefore there exists a
Radon measure $\mu$ such that $\Ge(\cdot,y)\weak\mu$ as
$\eps\to0$, weakly in the sense of measures. For $\varphi\in
C^\infty$ we compute:
\[
\varphi(y)=\eps^2\int\Ge(x,y)(-\Delta\varphi)(x)\dx+\int G(x,y)\varphi(x)\dx
\to\int\varphi\,\mathrm{d}\mu
\]
as $\eps\to0$. By density of $C^\infty$ in $C$, we conclude that
$\mu=\delta_y$. Proof of (ii). For $q=1$, we have:
\[
\|\Ge\ast h\|_1=\int|(\Ge\ast
h)(x)|\dx\le\int\dy|h(y)|\int\Ge(x,y)\dx=\int|h|=\|h\|_1.
\]
For $q=\infty$ we have, for any $x\in M$:
\[
|\Ge\ast
h(x)|\le\|h\|_\infty\int\Ge(x,y)\dy=\|h\|_\infty\int\Ge(x,y)\dx=\|h\|_\infty,
\]
and therefore $\|\Ge\ast h\|_\infty\le\|h\|_\infty$. The general
case follows by interpolation. Proof of (iii). 
Suppose $1<q<+\infty$. Let $\Ue=\Ge\ast
h$. Then we can write
\[
-\eps^2\Delta(\Ue-h)+(\Ue-h)=\eps^2\Delta h.
\]
Multiplying by $|\Ue-h|^{q-2}(\Ue-h)$ and integrating, we obtain
\begin{align*}
\eps^2(q-1)\int|\Ue-h|^{q-2}|\nabla(\Ue-h)|^2&+\int|\Ue-h|^q
\\
=&\eps^2\int\Delta h|\Ue-h|^{q-2}(\Ue-h).
\end{align*}
By positivity of the first term above and H\"older's inequality,
\[
\int|\Ue-h|^q\le\eps^2\int|\Delta h||\Ue-h|^{q-1}\le\eps^2\|\Delta
h\|_q\|\Ue-h\|_q^{q-1}.
\]
Hence $\|\Ue-h\|_q\le\eps^2\|\Delta h\|_q$ and (iii) follows
recalling the definition of $\Ue$ in the case $1<q<+\infty$.
Taking limits for $q\to1$ and $q\to+\infty$, we obtain the general case.
\end{proof}
\section*{Acknowledgments}
T.R.\ is grateful to Professor Gabriella Tarantello
for interesting and stimulating discussions.

\end{document}